\documentclass{amsart}

\usepackage{amsfonts,amssymb,latexsym}

\newtheorem{lemma}{Lemma}[section]

\newtheorem{theorem}[lemma]{Theorem}
\newtheorem{corollary}[lemma]{Corollary}

\begin{document}
\title{Sum-product estimates in finite fields via Kloosterman sums}
\author{Derrick Hart, \ Alex Iosevich, \ and \ Jozsef Solymosi}



\thanks{Research partially financed by NSERC and the National Science Foundation}

\begin{abstract} We establish improved sum-product bounds in finite fields using incidence theorems based on bounds for classical Kloosterman and related sums.
\end{abstract}

\maketitle

\tableofcontents

\section{Introduction}

\vskip.125in

Let $A \subset {\Bbb R}$. If $A$ is an arithmetic progression, then
$$ |A+A|=2|A|-1,$$ and
$$ |A \cdot A| \ge c{|A|}^{2-\epsilon},$$

where
$$ A+A=\{a+a': a,a' \in A\},$$
$$ A \cdot A=\{a \cdot a':a,a' \in A\},$$ and given a finite set $S$, $|S|$ denotes the number of its elements.

Similarly, if $B$ is an arithmetic progression and $A=\{2^n: n \in B\}$, then
$$ |A \cdot A|=2|A|-1,$$ and
$$ |A+A| \approx {|A|}^2.$$

Erd\H os and Szemer\'edi \cite{ESZ} proved the inequality

$$\max (|A+A|,|A\cdot A|)\geq c|A|^{1+\varepsilon}$$ for a small but positive $\varepsilon$, where $A$ is a subset of integers. They conjectured that
$$\max (|A+A|,|A\cdot A|)\geq c|A|^{2-\delta}$$ for any positive $\delta$.

After improvements in \cite{NA}, \cite{FO}, and \cite{CS} Elekes \cite{EL} showed that $\varepsilon\geq {1/4}$ if $A$ is a set of real numbers. His result was extended to complex numbers in \cite{CS} and \cite{SO}. For real and complex numbers the best known bound (\cite{Sol05}) says that 
$$ \max \{|A+A|, |A \cdot A|\} \ge c{|A|}^{\frac{14}{11}-\epsilon}.$$

See also, \cite{NT99}, \cite{C03} and \cite{ER03} for the discussion of the case when $|A+A|$ or $|A\cdot A|$ is small.

In the finite field setting the situation appears to be more complicated due to the fact that the Szemer\'edi-Trotter incidence theorem, the main tool in Euclidean setting, does not hold in the same generality and is, in general, much less well understood. It is known, however, via ground breaking work in \cite{BK03} and \cite{BKT04} that if $A \subset {\Bbb Z}_q$, $q$ a prime, than if $|A| \leq Cq^{1-\epsilon}$, for some $\epsilon>0$, then there exists $\delta>0$ such that
$$ \max \{|A+A|, |A \cdot A|\} \ge c{|A|}^{1+\delta}.$$

This bound does not yield a precise relationship between $\delta$ and $\epsilon$. The purpose of this paper is to establish a concrete value of $\delta$, in certain ranges of $|A|$, by using Kloosterman
sums, and to explore connections between this problem and that of incidences between points and circles in vector spaces over finite fields. We obtain reasonably good estimates when $|A|>>q^{\frac{1}{2}}$, and especially good ones when $|A| \approx q^{\frac{7}{10}}$. It would be great to obtain such estimates in the range $|A| \lesssim \sqrt{q}$, but this is out of our reach for the moment. See \cite{TV06} for a description of related results and applications to problems of additive combinatorics. Such estimates require assumptions on the existence of non-trivial subfields.  We note that in the range of exponents where the results of this paper are non-trivial, additional arithmetic assumptions are not required and it would be interesting to determine the precise parameters where this principle continues to hold.

Our main results are the following.
\begin{theorem} \label{main} Let $A \subset {\Bbb F}_q$, a finite field with $q$ elements. Suppose that
$$ |A+A|=m_1, \ |A \cdot A|=m_2.$$

Then
\begin{equation} \label{bigconclusion} {|A|}^3 \leq C(q^{-1}m_1^2m_2|A|+q^{\frac{1}{2}} m_1m_2).\end{equation}

In particular, if
$$ q^{\frac{1}{2}} \lesssim |A| \lesssim q^{\frac{7}{10}},$$ then
$$ \max \{|A+A|, |A \cdot A| \} \ge c\frac{|A|^{ \frac{3}{2} }}{q^\frac{1}{4}}.$$

\end{theorem}

Note that the best gain is achieved at the upper end of the range, when $|A| \approx q^{\frac{7}{10}}$. In this case
$$\max \{|A+A|, |A \cdot A| \} \ge c {|A|}^{\frac{8}{7}}.$$

Observe that when $|A| \approx q^{\frac{3}{4}},$ and the sumset is small, $|A+A|\leq C|A|,$ then the product set is large, $|A \cdot A| \ge cq.$

\begin{theorem} \label{main2} Let $A \subset {\Bbb F}_q$, a finite field with $q$ elements. Let $A^d=A \cdot A \dots \cdot A$, $d$ times. Suppose that
$$ |A+A|=m_1 \  \text{and} \ |A^d|=m_2.$$

Then
\begin{equation} \label{bigconclusion2} {|A|}^{2d} \leq C(q^{-1}{|A|}^d m_1^dm_2+q^{\frac{d-1}{2}}
{|A|}^{\frac{d}{2}} m_1^{\frac{d}{2}} m_2). 
\end{equation}

In particular, if $d=3$ and 
$$ q^{\frac{1}{2}} \lesssim |A| \lesssim q^{\frac{13}{21}},$$ then
$$ \max \{|A+A|, |A \cdot A \cdot A| \} \ge c\frac{{|A|}^{\frac{9}{5}}}{q^{\frac{2}{5}}}.$$

\end{theorem}

Note that the best gain is achieved at the upper end of the range, when $|A| \approx  q^{\frac{13}{21}}$. In that case 
$$\max \{|A+A|, |A \cdot A \cdot A| \} \ge c {|A|}^{\frac{15}{13}}.$$

This is a better exponent than the one given by Theorem \ref{main} in the context of triple products. Unfortunately, an application of the method of proof of Theorem \ref{main2} to products of higher degree does not yield reasonable results. We hope to address this problem in a subsequent paper. 

In analogy with our observation after the statement of Theorem \ref{main} above, we note the following consequence of Theorem \ref{main2}. 
\begin{corollary} \label{smallmember} Suppose that 
$$|A| \ge Cq^{\frac{1}{2}+\frac{1}{2d}}$$ and 
$$|A+A| \leq C'|A|.$$

Then 
$$|A^d| \ge cq.$$
\end{corollary} 

\vskip.125in

Our method also yields the following, related result.
\begin{theorem} \label{reacharound} Let $A \subset {\Bbb F}_q$, a finite field with $q$ elements. Suppose that $|A| \ge Cq^{\frac{1}{2}+\frac{1}{2d}}$ with a sufficiently large absolute constant $C$. Then
$$ (A-A) \cdot (A-A) \cdot \dots \cdot (A-A)={\Bbb F}_q,$$ where the product is taken $d$ times.
\end{theorem}

It is reasonable to conjecture that if $|A| \gtrsim q^{\frac{1}{2}+\epsilon}$, for some $\epsilon>0$, then
$(A-A) \cdot (A-A)={\Bbb F}_q$, but such a result is outside the realm of methods of this paper.

Moreover, it is easy to modify the proof of Theorem \ref{reacharound} so that $(A-A) \cdot \dots \cdot (A-A)$ may be replaced by $(A \pm A) \cdot \dots \cdot (A \pm A)$. We also recover the following corollary.

\begin{corollary} \label{subgroup} Let $H$ be a multiplicative subgroup of ${\Bbb F}_q^{*}$ such that
$|H| \ge Cq^{\frac{1}{2}+\epsilon}$, for some $\epsilon>0$. Then
$$ \pm kH={\Bbb F}_q^{*}$$ for $k \leq 2^{[\frac{1}{2 \epsilon}]}$, where $[t]$ denotes the smallest integer greater than $t$ and
$$ \pm kH=H \pm H \pm \dots \pm H,$$ $k$ times, with arbitrary choices of signs.
\end{corollary}

\vskip.125in

The main tool used in the proof of Theorem \ref{main} are properties of hyperbolas in vectors spaces over finite fields and incidence theorems (see Theorem \ref{circleincidence} below) based on Kloosterman sums bounds. Theorem \ref{main2} and Theorem \ref{reacharound} are based on higher dimensional incidence bounds given in Theorem \ref{multipleorgasms} below. Both are based on bounds for multi-dimensional Kloosterman sums obtained by Deligne.

\begin{theorem} \label{circleincidence} Let $E,F \subset {\Bbb F}^2_q$, the two dimensional vector space over a finite field with $q$ elements. Let $j \in {\Bbb F}_q^{*}$, the multiplicative group of ${\Bbb F}_q$. Then
$$ |\{(x,y) \in E \times F: (x_1-y_1)(x_2-y_2)=j\}| \leq C \left(q^{-1}|E||F|+\sqrt{q} \cdot \sqrt{|E||F|} \right).$$
\end{theorem}

The proof of analogous incidence theorems are contained in \cite{IR06} and \cite{IKo06} in the context of the Erd\H os distance problem in the vector spaces over finite fields. We shall give the argument below in the form required for the main result of this paper. The proof of the incidence bound is based on the following classical Kloosterman sum bound due to Andre Weil (\cite{W48}).

\begin{theorem} \label{kloosterman} Let ${\Bbb F}_q^{*}$ be as above and define
$$ K(a)=\sum_{t \in {\Bbb F}_q^{*}} \chi(at+t^{-1}),$$ where $\chi$ is a non-trivial additive character on 
${\Bbb F}_q$. Then if $a \not=0$,
\begin{equation} \label{kloostermanestimate} |K(a)| \leq 2 \sqrt{q}. \end{equation}

Moreover, the same estimate holds if the Kloosterman sum is "twisted" by a non-trivial multiplicative character. More precisely, (\ref{kloostermanestimate}) holds if $K$ is replaced by
$$ K_{\psi}(a)=\sum_{t \in {\Bbb F}_q^{*}} \chi(at+t^{-1}) \psi(t),$$ where $\psi$ is a non-trivial multiplicative character on ${\Bbb F}_q^{*}$.
\end{theorem}

\begin{theorem} \label{multipleorgasms} Let $E,F \subset {\Bbb F}_q^d$. Then 
$$ |\{(x,y) \in E \times F: (x_1-y_1)(x_2-y_2) \dots (x_d-y_d)=j \}|$$ 
$$ \leq C \left(|E||F|q^{-1}+q^{\frac{d-1}{2}} \sqrt{|E||F|} \right), $$ and if $j \not=0$ and $|E||F| \ge 
Cq^{d+1}$ with a sufficiently large constant $C$, then
$$ |\{(x,y) \in E \times F: (x_1-y_1)(x_2-y_2)\dots (x_d-y_d)=j \}|>0.$$
\end{theorem}

The main estimate in the proof of Theorem \ref{multipleorgasms} is the following result due to Pierre Deligne (\cite{D74}).
\begin{theorem} \label{multipleorgasmsest} Let $n \ge 1$, $\chi$ a non-trivial additive character, and define
$$ K_n(a)=\sum_{x \in {({\Bbb F}^{*}_n)}^n} \chi(x \cdot a+{(x_1x_2\dots x_n)}^{-1}).$$

Then if $a \not=(0, \dots, 0)$,
$$ |K_n(a)| \leq Cq^{\frac{n}{2}}.$$
\end{theorem}

\vskip.125in

\section{Proof of Theorem \ref{main}}

\vskip.125in

Let $A \subset {\Bbb F}_q$ and suppose that
$$ |A+A|=m_1,|A \cdot A|=m_2<<{|A|}^2=k.$$

By the pigeon-hole principle, there exists a large hyperbola in $A^2=A \times A$. More precisely, there exists $c \in {\Bbb F}_q^{*}$ such that
$$ |H_c| \ge \frac{k}{m_2},$$ where, without loss of generality, we may assume that 
$\frac{k}{m_2} \ge 5$, and
$$ H_c=\{(a, b) \in A^2: ab=c\} $$ which is a hyperbola.

\begin{lemma} \label{uniquehyperbola} If $|H_c|\geq 5$ then for any $(a_1,a_2) \in A \times A$,
$ (a_1,a_2)+H_c$ defines a unique hyperbola in ${\Bbb F}_q^2$.
\end{lemma}

To see this, suppose, for the sake of contradiction, that there exist $(a_1,a_2)$ and $(a'_1,a'_2)$, both in $A \times A$ such that $(a_1,a_2)+H_c$ and $(a'_1,a'_2)+H_c$ are both subsets of the same hyperbola $M$. Then for any $t \in M$, $(a_1-a'_1, a_2-a'_2)+t$ is a subset of $M$ as well as $(a_1-a'_1, a_2-a'_2)+it$ for any $i \in {\Bbb F}_q^{*}$. So, the hyperbola contains a line which is not possible.

Since $|A+A|=m_1$ it follows that we have $k$ hyperbolas containing at least $\frac{k}{m_2}$ points each on at most $m_1^2$ points. We are now ready to apply Theorem \ref{circleincidence} with $E= A \times A$ and $F=(A+A) \times (A+A)$. For every element $x \in E$ there are at least $\frac{k}{m_2}$ elements $y \in F$ such that ${(x_1-y_1)}{(x_2-y_2)}=c$. It follows that
$$ |E| \frac{k}{m_2} \leq C \left(q^{-1}|E| m_1^2+\sqrt{q} {|E|}^{\frac{1}{2}} m_1\right),$$ and (\ref{bigconclusion}) follows. This completes the proof of Theorem \ref{main}.

\vskip.125in

\section{Proof of Theorem \ref{main2}}

\vskip.125in 

By the pigeon-hole principle, there exists a large hyperbola in the $d$-fold product $A^d$. More precisely, there exists $c \in {\Bbb F}_q^{*}$ such that
$$ |H_c| \ge \frac{k}{m_2},$$ where, without loss of generality, we may assume that 
$\frac{k}{m_2} \ge d+3$, and
$$ H_c=\{(a_1, a_2, \dots, a_d) \in A^d: a_1a_2\dots a_d=c\}. $$ 

\begin{lemma} \label{uniquehyperbola3} If $|H_c| \geq d+3$ then for any $x \in A \times A \times \dots \times A$, $x+H_c$ defines a unique hyperbola in ${\Bbb F}_q^d$. 
\end{lemma}

To see this, suppose, for the sake of contradiction, that there exist $(a_1,a_2, \dots, a_d)$ and 
$(a'_1,a'_2, \dots, a'_d)$, both in $A \times A \times \dots \times A$ such that $(a_1,a_2, \dots, a_d)+H_c$ and $(a'_1,a'_2, \dots, a'_d)+H_c$ are both subsets of the same hyperbola $M$. Then for any $m \in M$, $(a_1-a'_1, a_2-a'_2, \dots, a_d-a'_d)+m$ is a subset of $M$ as well as $(a_1-a'_1, a_2-a'_2, \dots a_d-a'_d)+im$ for any $i \in {\Bbb F}_q^{*}$. So, the hyperbola contains a line which is not possible. 

Since $|A+A|=m_1$ it follows that we have $k$ hyperbolas containing at least $\frac{k}{m_2}$ points each on at most $m_1^d$ points. We are now ready to apply Theorem \ref{multipleorgasms} with
$E=A^d$ and $F=(A+A) \times (A+A) \times \dots \times (A+A)$. For every element $x \in E$ there are at least $\frac{k}{m_2}$ elements $y \in F$ such that $(x_1-y_1)(x_2-y_2) \dots (x_d-y_d)=c$. It follows that
$$ |E|\frac{k}{m_2} \leq C \left(q^{-1}|E| m_1^d+q^{\frac{d-1}{2}} {|E|}^{\frac{1}{2}} 
m_1^{\frac{d}{2}}\right),$$ and (\ref{bigconclusion2}) follows. This completes the proof of Theorem 
\ref{main2}.

\vskip.125in

\section{Proof of Theorem \ref{circleincidence}}

\vskip.125in

We have
$$ |\{(x,y) \in E \times F: (x_1-y_1)(x_2-y_2)=j\}|$$
\begin{equation} \label{kloostermansetup}=\sum_{x,y \in {\Bbb F}_q^2} E(x)F(y)S_j(x-y),\end{equation} where $E$ and $F$ are characteristic functions of $E$ and $F$ respectively and $S_j$ is the characteristic function of the set
$$ \{x \in {\Bbb F}_q^2: x_1x_2=j\}.$$

Recall (see e.g. \cite{IK04} that if $f$ is a function from ${\Bbb F}_q^2$ to the complex numbers, then
$$ \widehat{f}(m)=q^{-2} \sum_{x \in {\Bbb F}_q^2} \chi(-x \cdot m) f(x),$$ the Fourier transform of $f$, where $\chi$ is a non-trivial additive character on ${\Bbb F}_q$. Also recall that
$$ f(x)=\sum_{m \in {\Bbb F}_q^2} \chi(x \cdot m) \widehat{f}(m),$$ and
\begin{equation} \label{plancherel} \sum_{m \in {\Bbb F}_q^2} {|\widehat{f}(m)|}^2=q^{-2} \sum_{x \in {\Bbb F}_q^2} {|f(x)|}^2. \end{equation}

It follows that the expression in (\ref{kloostermansetup}) equals
$$ \sum_{x,y,m \in {\Bbb F}_q^2} \chi((x-y) \cdot m) E(x)F(y) \widehat{S}_j(m)$$
$$=q^4 \sum_{m \in {\Bbb F}_q^2} \widehat{E}(m) \cdot \overline{\widehat{F}}(m) 
\cdot \widehat{S}_j(m)$$
$$=q^4 \cdot q^{-2} \cdot \widehat{E}(0,0) \cdot q^{-2} \widehat{F}(0,0) \cdot q^{-2} \widehat{S}_j(0,0)$$
$$ +q^4 \sum_{m \not=(0,0)} \widehat{E}(m) \cdot \overline{\widehat{F}}(m) \cdot \widehat{S}_j(m)$$
$$=q^{-2} \cdot |E| \cdot |F| \cdot |S_j|+q^4 \sum_{m \not=(0,0)} \widehat{E}(m) \cdot \overline{\widehat{F}}(m) \cdot \widehat{S}_j(m)=I+II.$$

We shall need the following result which we prove at the end of the section.
\begin{lemma} \label{sphere} Suppose that $j \not=0$. With the notation above,
$$ \# S_j=q+O(\sqrt{q}),$$ and
\begin{equation} \label{fouriersphere} |\widehat{S}_j(m)| \leq Cq^{-\frac{3}{2}}\end{equation} provided that $m \not=(0,0)$.
\end{lemma}

We have
$$ I \leq C|E||F|q^{-1},$$ whereas the application of Cauchy-Schwarts shows that
\begin{equation} \label{cs} II \leq q^4 {\left( \sum_{m \in {\Bbb F}_q^2}
{|\widehat{E}(m)|}^2 \right)}^{\frac{1}{2}} \cdot {\left( \sum_{m \in {\Bbb F}_q^2}
{|\widehat{F}(m)|}^2 \right)}^{\frac{1}{2}} \cdot \sup_{m \not=(0,0)} |\widehat{S}_j(m)|. \end{equation}

Applying (\ref{plancherel}) to the first two terms in (\ref{cs}) and using (\ref{fouriersphere})
of Lemma \ref{sphere} to estimate the third term, we see that
$$ II \leq C q^4 \cdot q^{-2} \cdot \sqrt{|E||F|} \cdot q^{-\frac{3}{2}}=Cq^{\frac{1}{2}} \cdot {|E|}^{\frac{1}{2}} \cdot {|F|}^{\frac{1}{2}},$$
and the proof of Theorem \ref{circleincidence} is complete up to the proof of Lemma \ref{sphere}.

\vskip.125in

\subsection{Proof of Lemma \ref{sphere}} We have
$$ |S_j|=\sum_{x \in {\Bbb F}_q^2} S_j(x)=\sum_{x \in {\Bbb F}_q^2} q^{-1} \sum_{t \in {\Bbb F}_q} \chi(t(x_1x_2-j))$$
$$=q+q^{-1} \sum_{t \not=0} \chi(-tj) \sum_{x \in {\Bbb F}_q^2} \chi(t(x_1x_2))=q+D(q).$$

Let $x_1=\frac{y_1+y_2}{2}$ and $x_2=\frac{y_1-y_2}{2}$. We then see that
$$ \sum_{x \in {\Bbb F}_q^2} \chi(t(x_1x_2))=\sum_{y \in {\Bbb F}_q^2} \chi(t(y_1^2-y_2^2))=q \cdot \psi_q(t),$$ where $\psi_q$ is a non-trivial multiplicative character of ${\Bbb F}_q^{*}$. See, for example, \cite{IK04}.

It follows that
$$ |D(q)| \leq \left| \sum_{t \not=0} \chi(-tj) \psi_q(t) \right| \leq C \sqrt{q}, $$ by a classical theorem on Fourier transforms of non-trivial multiplicative characters. See e.g. \cite{IK04}. This completes the proof of the first part of Lemma \ref{sphere}.

To prove (\ref{fouriersphere}) we write
$$ \widehat{S}_j(m)=q^{-2} \sum_{x \in {\Bbb F}_q^2} q^{-1} \sum_{t \in {\Bbb F}_q} \chi(-x \cdot m) \chi(t(x_1x_2-j))$$
\begin{equation} \label{nearejaculation}=q^{-3} \sum_{x \in {\Bbb F}_q^2} \sum_{t \not=0} \chi(-x \cdot m) \chi(t(x_1x_2-j)). \end{equation}

Changing variables as above and completing the square, we see, as before, that
$$ \sum_{x \in {\Bbb F}_q^2} \chi(-x \cdot m) \chi(t(x_1x_2))$$
$$= \sum_{y \in {\Bbb F}_q^2} \chi \left( -\frac{y_1+y_2}{2}m_1-\frac{y_1-y_2}{2}m_2 \right) \chi(t(y_1^2-y_2^2))$$
$$=q \chi \left( \frac{m_1^2-m_2^2}{4t} \right) \psi_q(t).$$

Plugging this into (\ref{nearejaculation}) we get
$$ q^{-2} \sum_{t \not=0} \chi \left( \frac{m_1^2-m_2^2}{4t}+tj \right) \psi_q(t),$$ and the conclusion follows from Theorem \ref{kloosterman} provided that $m_1^2-m_2^2 \not=0$. If $m_1^2-m_2^2=0$, the desired estimate follows from the aforementioned result on Fourier transforms of non-trivial multiplicative characters. The proof of Lemma \ref{sphere} is complete.

\vskip.125in

\section{Proof of Theorem \ref{reacharound}}

\vskip.125in

Suppose for a moment that Theorem \ref{multipleorgasms} holds. Let $E=A^d$. It then follows immediately that for any $j \in {\Bbb F}_q^{*}$, there exists $x, y \in A^d$ such that
$$ (x_1-y_1)(x_2-y_2) \dots (x_d-y_d)=j,$$ which means that
$$ (A-A) \cdot (A-A) \dots (A-A)={\Bbb F}_q, $$ as desired.

This matters are reduced to proving Theorem \ref{multipleorgasms} which is where we now turn our attention. Let
$$ M_j=\{x \in {\Bbb F}_q^d: x_1x_2 \dots x_d=j\}.$$

We have, following the proof of Theorem \ref{circleincidence},
$$ |\{(x,y) \in E \times F: (x_1-y_1) \dots (x_d-y_d)=j\}|$$
$$=\sum_{x,y \in {\Bbb F}_q^d} E(x)F(y)M_j(x-y)$$
$$=q^{-d}|E||F| \cdot |M_j|+q^{2d} \sum_{m \not=(0, \dots, 0)} \widehat{E}(m) \overline{\widehat{F}(m)} \widehat{M}_j(m)=I+II.$$

\begin{lemma} \label{orgasminaction} With the notation above,
\begin{equation} \label{countM} |M_j|=q^{d-1}+O(q^{d-2}), \end{equation} and if $m \not=(0, \dots, 0)$,
\begin{equation} \label{fourierM} |\widehat{M}_j(m)| \leq Cq^{-\frac{d+1}{2}}. \end{equation}
\end{lemma}

We postpone the proof of the lemma for a moment and complete the proof of Theorem \ref{multipleorgasms}. Using (\ref{countM}) we see that
$$C_1 q^{-1} |E||F| \leq I \leq C_2 q^{-1} |E||F|. $$

Using (\ref{fourierM}) and Cauchy-Schwartz we see that
$$ |II| \leq C_3 q^{2d} \cdot q^{-d} \cdot \sqrt{|E||F|} \cdot q^{-\frac{d+1}{2}}=
q^{\frac{d-1}{2}} \cdot \sqrt{|E||F|}. $$

This establishes the first part of Theorem \ref{multipleorgasms}. Since we have a lower bound on $I$ and an upper bound on $II$, we conclude that if $|E| \ge Cq^{\frac{d+1}{2}}$ with $C$ sufficiently large, then $I+II>0$ and the second part of Theorem \ref{multipleorgasms} follows.

\vskip.125in

\subsection{Proof of Lemma \ref{orgasminaction}} Let $\chi$ be a non-trivial additive character. We have
$$ |M_j|=\sum_{x_1 \dots x_d=j} 1=\sum_{x \in {\Bbb F}_q^d} q^{-1} \sum_{t \in {\Bbb F}_q} \chi(t(x_1 \dots x_d-j))$$
$$=q^{d-1}+q^{-1} \sum_{t \not=0} \chi(-tj) \sum_{x \in {\Bbb F}_q^d} \chi(tx_1 \dots x_d)=q^{d-1}+O(q^{d-2})$$ thus establishing (\ref{countM}).

We now prove (\ref{fourierM}). We have
$$ \widehat{M}_j(m)=q^{-d} \sum_{x_1 \dots x_d=j} \chi(-x \cdot m)$$
$$=q^{-d} \sum_{x_1, x_2 \dots, x_{d-1}} \chi(-x_1m_1-\dots-x_{d-1}m_{d-1}-m_d{(x_1 \dots x_{d-1})}^{-1}),$$ and the conclusion follows from Theorem \ref{multipleorgasmsest}. This completes the proof of Lemma \ref{orgasminaction}.


\end{document}